# Optimization of the marketing mix in the health care industry

Dominique Haughton, Bentley University, Université Paris I, and Université Toulouse I[1]
Guangying Hua, Boston Consulting Group
Danny Jin, Epsilon
John Lin, Epsilon
Qizhi Wei, Epsilon
Changan Zhang, Bentley University[2]

Abstract. This paper proposes data mining techniques to model the return on investment from various types of promotional spending to market a drug and then uses the model to draw conclusions on how the pharmaceutical industry might go about allocating marketing expenditures in a more efficient manner, potentially reducing costs to the consumer.

---

[1] Supported by NSF DMS-1106388
[2] Supported by NSF DMS-1106388

# Optimization of the marketing mix in the health care industry

This paper proposes data mining techniques to model the return on investment from various types of promotional spending to market a drug and then uses the model to draw conclusions on how the pharmaceutical industry might go about allocating marketing expenditures in a more efficient manner, potentially reducing costs to the consumer. The main contributions of the paper are twofold. First, it demonstrates how to undertake a media mix optimization process in the pharmaceutical context and carry it through from the beginning to the end. Second, the use of Directed Acyclic Graphs (DAGs) to help unravel the effects of various marketing media on sales volume, notably direct and indirect effects, is proposed.

## Introduction; current best practice in the pharmaceutical industry

The matter of marketing expenditures by the pharmaceutical industry continues to attract considerable attention (see editorial by Mukherjee, 2012 and the introduction by Crié and Chebat to the special issue in the Journal of Business Research on health marketing, 2013), particularly at a time when health care costs are coming under tightened scrutiny. In this context, while marketing expenditures in the pharmaceutical industry have continued to decline (Eyeforpharma, 2014), concerns remain about the proportion of such expenditures relative to that of Research and Development (Pew 2013), and about direct-to-consumer (DTC) marketing programs, legal only in the United States and New Zealand (Liang and Mackey 2011).

Marketing efforts in the pharmaceutical industry can be grouped into two categories: 1) direct-to-consumer (DTC) advertising such as TV, print, internet ads as well as Customer Relationship Management (CRM) programs, and 2) professional or direct-to-physician (DTP) promotions including samples, details, journal advertising, professional interactions, and non-personal promotions (NPP). Most of the market mix analyses in the industry focus on measuring the effects of the different DTC advertising channels after controlling for DTP promotions. These analyses are usually performed by creating a control group of physicians in Designated Market Areas *(*DMAs*,* geographic areas defined by the Nielsen Media Research Company as a group of counties that make up a particular television market) or other geographic areas matched to a test group with identical levels of DTP promotions and physician characteristics. The test and control groups have varying levels of DTC activities; this allows for the measurement of the effect of the different DTC efforts.

New product launches typically utilize market mix models to simulate potential outcomes (market volume, ROI, etc.) under different media spending scenarios. These models incorporate intermediate outcomes such as brand awareness and other consumer behavior that result from consumers asking their doctor about the brand and subsequently filling prescription for this brand.

We next review pertinent recent literature and outline the objectives and contributions of this paper.

**Review of recent literature**

While the issue of optimizing marketing mixes has been of interest for many years in the marketing literature, relatively little such published work is dedicated to the health care industry. Of particular relevance to this paper, albeit in another industry (tourism), is the work by Wolfe and Crotts (2011) where a marketing mix methodology is applied to ticket sales for a museum.

In general the literature in health care marketing makes the point that decisions on pharmaceutical marketing expenditures are often of a qualitative or policy nature, but that it is in fact important to quantitatively link marketing expenditures and new prescription volume, since a more optimal allocation of expenditures would lead to savings that could be expected to decrease the price of drugs for patients.

Past work has attempted to address the issue of correlated predictors and the problem of a suitable choice of a set of predictors for the number of new prescriptions for a drug, or any other yield dependent variable. Lim and Kirikoshi (2008) and Lim, Kirikoshi and Okano (2008) propose to use a genetic algorithm (with the predictive power of a model as an objective function) to help select a suitable set of predictors, in combination with a neural net, or a partial least squares regression. The models perform well, but one issue that remains to be addressed is of how to operationalize the model for marketing expenditure allocation, when, for example, some coefficients are negative. In many cases, practitioners faced with inconvenient negative coefficients are forced to ignore them. However, negative coefficients can and sometimes do represent reality; we will return to this issue later in the paper.

Manchanda et al. (2008) propose a model which takes into account the fact that the level of marketing-mix variables is often set with at least partial a-priori knowledge of the likely level of response for each variable, and find that this approach results in a more precise estimation of response parameters. They also find that physicians are often not detailed optimally, but that high-volume physicians are detailed to a greater extent than low-volume physicians without regard to responsiveness. The work by Singh (2008) proposes a conceptual model that incorporates several aspects of the network created by physician-salesperson dyads.

More recently, Huber et al. (2012) have constructed models, applicable to different categories of Over The Counter (OTC) products, which help to identify drivers of product sales in the OTC market. Gönül and Carter (2012) construct a model of physician prescribing volume of older versus newer drugs and confirm the importance of sampling and detailing.

It is useful to note, since markets do differ across countries, that the context discussed in this paper is that of pharmaceutical marketing in the United States. In particular, from a research point of view, data availability varies across countries. Obtaining data for research about marketing pharmaceutical expenditures is notoriously difficult in the United States, and can occur almost exclusively in the context of collaborations between academics and practitioners and by of course masking all identifiers of drugs, patients, physicians and even pharmaceutical organizations involved.

In summary, past literature has identified main drivers of prescribing behavior in a number of pharmaceutical industry cases, and preliminary marketing mix models have been constructed. However, the issue of correlated predictors remains a challenge, as well as the operationalization of media mix models in the pharmaceutical context. In order to address this gap, this paper makes the following twofold contributions:

1. For the first time in the literature (to the best of our knowledge) we demonstrate how to undertake a media mix optimization process and carry it through from the beginning to the end; to achieve this goal, synthetic data are employed.
2. We propose the use of Directed Acyclic Graphs (DAGs) to help unravel the effects of various marketing media on sales volume, notably direct and indirect effects.

We now outline the two main phases which are typically included in a media mix analysis and, after providing details about the synthetic dataset used here, we describe and implement each step of the two phases in a prototype marketing mix process.

**Modeling approach**

The approach to optimizing the marketing mix for a product typically proceeds in the following two phases:

**Phase 1**

A model is built for the output variable (typically new prescriptions in a given time period) in terms of a number of relevant independent variables. Phase 1 corresponds to the central and upper sections of Figure 1.

**Phase 2**

Once a satisfactory model has been constructed, the model is employed to evaluate the contribution of each marketing activity to the new prescriptions. The objective is to obtain output such as obtained in Figure 2, although the marketing components will tend to differ in a pharmaceutical context. This step corresponds to the central lower part of Figure 1.

*Figure 1: Framework for Marketing ROI assessment and optimization*

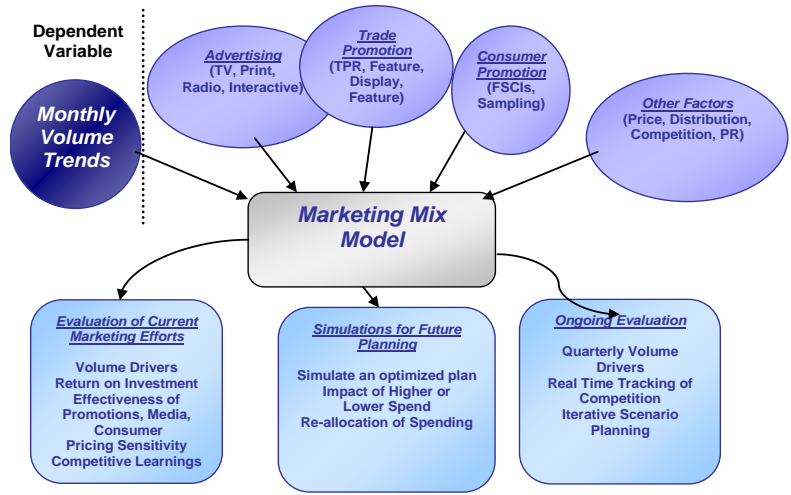

*Figure 2: Typical output from a marketing mix modeling analysis (extracted from Wikipedia, Marketing Mix Modeling, 2014)*

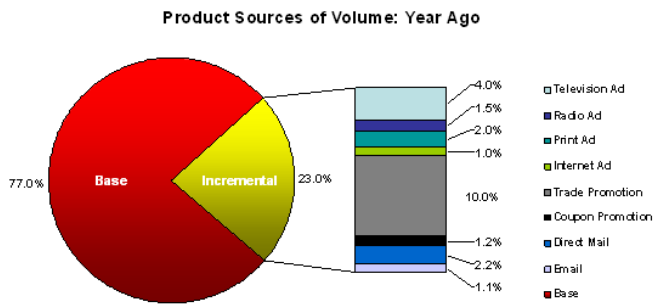

In order to illustrate the methods proposed in this paper, we describe below a "Prototype modeling approach for Phases 1 and 2" on the basis of synthetic data. We next introduce the data used to that effect.

## Data

Using the Stata statistical software, we generate synthetic data simulated from the correlation matrix and summary statistics for 12 variables in Lim et al. (2008). The variables include, over a period of 71 months, the number of new prescriptions (*nrx*) and eleven variables related to various aspects of different advertising activities (see Table 1). In Lim et al. (2008) the dataset was obtained from a marketing consultant and concerns an antibiotic drug in the United States.

Complete information about the variables can be found in Lim et al. (2008), but we reproduce the definitions given in that paper here for the sake of completeness and convenience. As defined by Lim et al. (2008):

> *"Contacts" (CON) is a product-level report of promotional actions that is provided to physicians and "calls" (CAL) measures the total number of visits made by pharmaceutical representatives to physicians. A CON can be a full product discussion with a physician, a drug fair set up at the hospital for physicians or a delivery of a product sample. Several products may be discussed during a single call, resulting in the possibility of multiple CONs in a CAL. "Cost of contacts" (COC) includes the costs associated with detailing of representatives that are directed to physicians. "Cost per contact" (CPC) is an estimate of cost per contact whereas "minutes" (MIN) is the projected sum of time spent with physicians. In a broad sense we interpret CPC as the quality of CON. This is a rough approximation that accounts for the difficulty and complexity of assessing the physician's overall impression of the representative's visits. "Journal advertising spending" (JAS) captures the expenditure of advertising in medical journals. "Ads" (ADS) measures the number of different layouts of product advertisements in medical journals. If the same ad appears in two journals, it is counted twice. "Ad pages circulated" (ADP) represents the number of total ad pages circulated in journals for a particular product. "Sample" (SAM) shows the projected volume of a product provided as samples to physicians whereas "extended units samples" (EUS) measures the amount of a product sampled as the number of packages multiplied by the size of the package in tablets, capsules, milliliters, etc. EUS is appropriate for use when the products being compared are similar in terms of dosage form. "Retail value of sample" (RVS) represents the retail value of SAM. "New prescription volume" (NRx) represents the count of new prescriptions dispensed by pharmacists.*

We note that in addition to data on the brand of antibiotic (referred to as A) simulated from the statistical summaries in Lim et al. (2008), it would be potentially possible to simulate data from the statistical summaries on three other brands B, D and D also covered in the work by Lim et al. (2008).

**Analysis: phase 1**

As is readily visible on Figure 3, a major issue when attempting to build a model for the number of new prescriptions in terms of the eleven predictors is that a high level of correlation exists among the predictors. For example, contacts and cost of contacts are highly correlated as one might expect. Contacts, calls and cost of contacts are highly correlated to the retail value of samples.

The first row of the matrix scatter plot graphs the dependent variable *nrx* in terms of the eleven predictors; a positive linear relationship emerges between *nrx* and each predictor except for *cpc* (Cost per Contact) where that relationship is negative.

*Table 1: List of synthetic variables (following Lim, Kirikoshi and Okano (2008)) and summary statistics*

|  | N | Minimum | Maximum | Mean | Std. Deviation |
|---|---|---|---|---|---|
| Con (Contacts) | 71 | 7657.00 | 99188.00 | 54113.88 | 21606.72 |
| Cal (Calls) | 71 | 5361.00 | 77657.00 | 38964.52 | 17162.88 |
| Coc (Cost of Contacts) | 71 | 737297.00 | 8041793.00 | 4357317.52 | 1740648.34 |
| Cpc (Cost per Contact) | 71 | 62.00 | 114.00 | 87.92 | 10.67 |
| Min (Minutes) | 71 | 32682.00 | 334945.00 | 182602.33 | 72862.46 |
| Jas (Journal Advertising Spending) | 71 | .00 | 486943.00 | 217926.80 | 112177.79 |
| Ads (Ads) | 71 | .00 | 30.00 | 14.12 | 6.57 |
| Adp (Ad Pages Circulated) | 71 | .00 | 82.00 | 39.16 | 20.22 |
| Sam (Samples) | 71 | 615.00 | 2614948.00 | 1094355.12 | 627396.75 |
| Eus (Extended Unit Samples) | 71 | 1230.00 | 8204439.00 | 3556610.79 | 1743431.60 |
| Rvs (Retail Value of Sample) | 71 | 6660.00 | 18502700.00 | 8582259.60 | 3820827.20 |
| Nrx (New Prescription Volume) | 71 | 84895.00 | 2466777.00 | 1080544.09 | 512578.60 |

If we attempt to run a linear regression for *nrx* in terms of the eleven predictors, several variables fail to be significant and the negative coefficients for some of the variables might be considered counter-intuitive (Table 2). A stepwise regression for *nrx* (Table 3) yields a negative coefficient for *jas* and includes a small number of significant predictors.

One option would be to remove *jas* and evaluate the resulting model; this yields an R-square of .910, not a very large loss from the value .919, and the coefficients of the remaining three variables are positive and significant. However, and this is a matter that would need to be elucidated from domain knowledge of the data, the negative sign could in fact represent a reality. We are aware of an interesting case where a negative sign appeared in a similar marketing mix, but this time in the banking industry. It transpired that all coefficients were positive in what looked like a sensible linear regression model, except that of newspaper advertising expenditures. A more careful examination revealed that interest rates were higher for the bank under consideration than for competitors and that those rates were published in the same newspapers!

*Table 2: Linear regression of nrx in terms of all eleven predictors (R-square=.931)*

| Model | | Unstandardized Coefficients | | Standardized Coefficients | | |
|---|---|---|---|---|---|---|
| | | B | Std. Error | Beta | t | Sig. |
| 1 | (Constant) | 1296725.732 | 522210.488 |  | 2.483 | .016 |
| | con | -31.255 | 14.707 | -1.318 | -2.125 | .038 |
| | cal | 12.363 | 16.948 | .414 | .729 | .469 |
| | coc | .190 | .119 | .646 | 1.595 | .116 |
| | cpc | -11301.270 | 4710.481 | -.235 | -2.399 | .020 |
| | min | 1.228 | 1.648 | .174 | .745 | .459 |
| | jas | -1.556 | .583 | -.341 | -2.668 | .010 |
| | ads | 33979.885 | 14621.334 | .435 | 2.324 | .024 |
| | adp | 2713.728 | 5320.917 | .107 | .510 | .612 |
| | sam | .480 | .239 | .588 | 2.009 | .049 |
| | eus | -.137 | .077 | -.465 | -1.782 | .080 |
| | rvs | .075 | .025 | .563 | 3.034 | .004 |

*Figure 3: Matrix scatter plot for twelve synthetic variables*

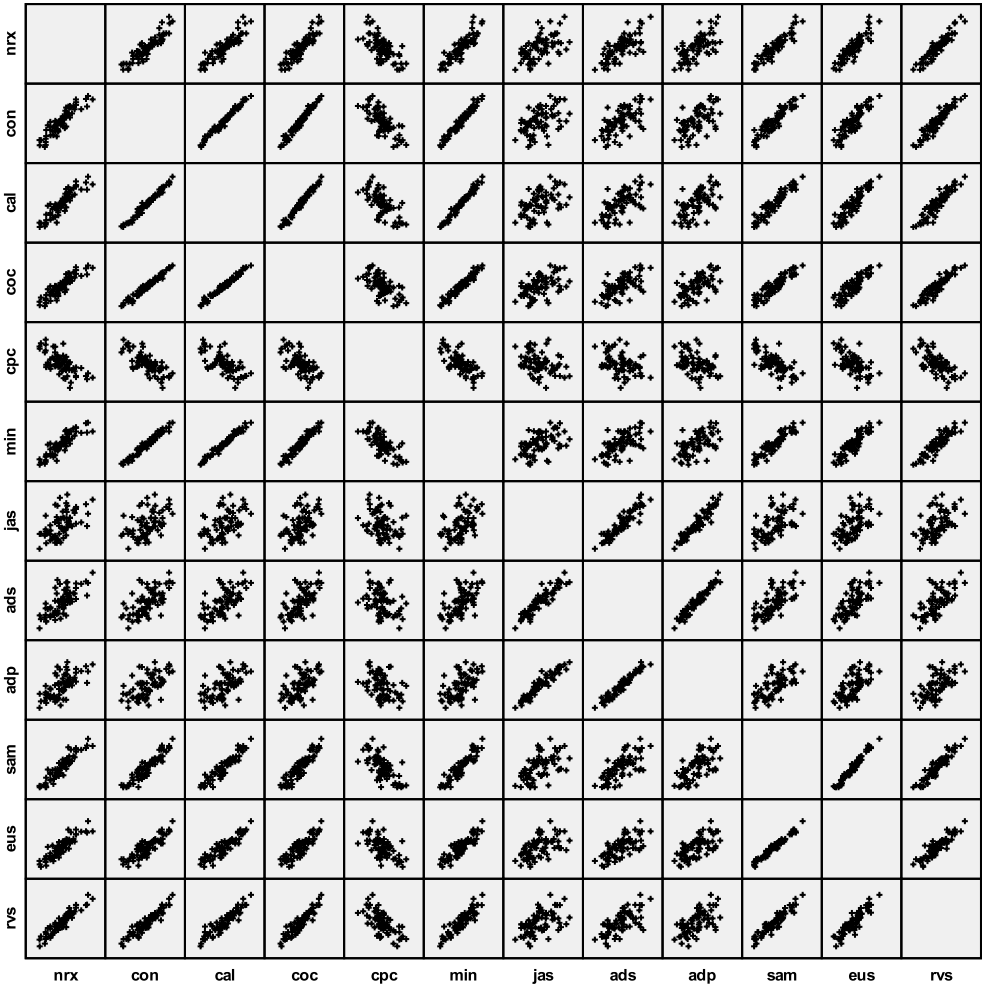

*Table 3: Stepwise linear regression results for nrx (R-square=.919)*

| Model | | Unstandardized Coefficients | | Standardized Coefficients | | |
|---|---|---|---|---|---|---|
| | | B | Std. Error | Beta | t | Sig. |
| 1 | (Constant) | 22055.430 | 52669.622 | | .419 | .677 |
| | jas | -1.214 | .450 | -.266 | -2.699 | .009 |
| | ads | 33026.816 | 9082.479 | .423 | 3.636 | .001 |
| | sam | .228 | .088 | .279 | 2.602 | .011 |
| | rvs | .071 | .014 | .527 | 5.036 | .000 |

A potentially better approach is to extract uncorrelated factors to explain as large as possible an amount of variability in the eleven predictors. A solution with two factors, rotated with a quartimax rotation, explains 91.3% of the total variability in the predictors and is presented in Table 4.

*Table 4: Factor analysis of all eleven predictors*

| Component | Initial Eigenvalues | | | Extraction Sums of Squared Loadings | | | Rotation Sums of Squared Loadings | | |
|---|---|---|---|---|---|---|---|---|---|
| | Total | % of Variance | Cumulative % | Total | % of Variance | Cumulative % | Total | % of Variance | Cumulative % |
| 1 | 8.487 | 77.153 | 77.153 | 8.487 | 77.153 | 77.153 | 7.938 | 72.160 | 72.160 |
| 2 | 1.611 | 14.649 | 91.801 | 1.611 | 14.649 | 91.801 | 2.106 | 19.148 | 91.308 |
| 3 | .476 | 4.325 | 96.126 | .476 | 4.325 | 96.126 | .530 | 4.819 | 96.126 |
| 4 | .220 | 2.001 | 98.127 | | | | | | |
| 5 | .097 | .885 | 99.012 | | | | | | |
| 6 | .051 | .467 | 99.479 | | | | | | |
| 7 | .022 | .198 | 99.677 | | | | | | |
| 8 | .021 | .192 | 99.868 | | | | | | |
| 9 | .009 | .083 | 99.951 | | | | | | |
| 10 | .003 | .031 | 99.982 | | | | | | |
| 11 | .002 | .018 | 100.000 | | | | | | |

*Table 5: Loadings of each predictor onto each rotated factor*

| | Component | |
|---|---|---|
| | 1 | 2 |
| con | .988 | .046 |
| cal | .984 | .053 |
| coc | .964 | .066 |
| cpc | -.768 | .098 |
| min | .979 | .031 |
| jas | .437 | .873 |
| ads | .594 | .787 |
| adp | .539 | .828 |
| sam | .958 | .169 |
| eus | .921 | .193 |
| rvs | .951 | .137 |

Extraction Method: Principal Component Analysis.
Rotation Method: Quartimax with Kaiser Normalization.
a. Rotation converged in 3 iterations.

Table 5 reveals that the first factor seems to represent activities directly targeted to physicians, while the second factor represents more general advertising activities. Table 6 displays the coefficients for the Z score of each predictor needed to compute the factors, and Table 7 presents the results of a regression of *nrx* on the two factors (R-square = .885).

One can then extract a formula for the estimated new prescription volume in terms of the Z scores of the predictors (and ultimately in terms of the predictors themselves, of course) and use an optimization process, here a simple linear programming application, to find the optimal value of the Z scores under some constraints, such as total budget, etc. For illustrative purposes, we set the lower limit at -2 and the upper limit at 4 for each Z score, with an upper limit of 30 for the sum of all Z scores, and obtain the results in Table 8.

*Table 6: Coefficients to compute each factor in terms of the Z-scores of the predictors*

|     | Component 1 | Component 2 |
|-----|-------------|-------------|
| con | .151        | -.109       |
| cal | .149        | -.105       |
| coc | .144        | -.094       |
| cpc | -.135       | .162        |
| min | .151        | -.117       |
| jas | -.051       | .446        |
| ads | -.014       | .375        |
| adp | -.028       | .406        |
| sam | .129        | -.034       |
| eus | .120        | -.015       |
| rvs | .132        | -.051       |

Extraction Method: Principal Component Analysis.
Rotation Method: Quartimax with Kaiser Normalization.
Component Scores.

*Table 7: Linear regression of nrx in terms of the two factors*

**Model Summary**

| Model | R      | R Square | Adjusted R Square | Std. Error of the Estimate |
|-------|--------|----------|-------------------|----------------------------|
| 1     | .941[a]| .885     | .882              | 1.76432E5                  |

**Coefficients[a]**

| Model |            | Unstandardized Coefficients | | Standardized Coefficients | t      | Sig. |
|-------|------------|-----------|------------|------------|--------|------|
|       |            | B         | Std. Error | Beta       |        |      |
| 1     | (Constant) | 1080544.093 | 20938.592 |            | 51.605 | .000 |
|       | Factor1    | 470291.077 | 21087.623 | .918       | 22.302 | .000 |
|       | Factor2    | 106415.943 | 21087.623 | .208       | 5.046  | .000 |

a. Dependent Variable: nrx

Two remarks are in order. First, while we have used a factor analysis on the predictors, and then used the factors in a regression analysis for new prescriptions, we could also use a Partial Least Squares regression, which will attempt to extract factors which explain a satisfactory proportion of the variability in the predictors as well as in the new prescriptions. Such an analysis performed on the synthetic data yields similar results to those above, but with a slightly smaller of proportion of the variability in the predictors explained (90.5%) and a bit higher a proportion of variability in new prescription explained (89.3%). Second, our analysis does not take into account any correlation in the errors of the regression model. The synthetic data were created without such correlation, but it is very likely that error correlation exists in real data, and would need to be taken into account, for example by modeling the errors as autoregressive or even ARMA (Auto Regressive Moving Average).

## Analysis: phase 2

Table 8 reveals that an optimal mix would set the variables *con, cal, coc, min, ads, sam, eus* and *rvs* at their maximum levels, the variables *cpc* and *jas* at their minimum level, and the variable *adp* at a high but not maximum level (2 standard deviations above its mean).

*Table 8: Optimization of predictors using linear programming*

| Variable | Coefficient of Z Score of variable in objective function | Optimal Z Score | Contribution of variable to objective function | Lower limit | Upper limit |
|---|---|---|---|---|---|
| con | 59167.95 | 4 | 236671.81 | -2 | 4 |
| cal | 58912.58 | 4 | 235650.32 | -2 | 4 |
| coc | 57713.02 | 4 | 230852.08 | -2 | 4 |
| cpc | -46433.24 | -2 | 92866.49 | -2 | 4 |
| min | 58693.49 | 4 | 234773.95 | -2 | 4 |
| jas | 23596.74 | -2 | -47193.47 | -2 | 4 |
| ads | 33295.52 | 4 | 133182.09 | -2 | 4 |
| adp | 29855.06 | 2 | 59710.12 | -2 | 4 |
| sam | 56987.67 | 4 | 227950.67 | -2 | 4 |
| eus | 54746.69 | 4 | 218986.76 | -2 | 4 |
| rvs | 56685.92 | 4 | 226743.68 | -2 | 4 |
|  |  |  |  |  |  |
|  |  |  | **Objective function** |  |  |
|  | **Constraint on the sum of Z Scores for all variables** | 30 | 1850194.49 |  |  |

## Another approach to unraveling the effects of correlated predictors: Directed Acyclic Graphs

Another approach to help unravel the effects of correlated predictors on the new prescription volume is the construction of a Directed Acyclic Graph (DAG), followed by the estimation of a structural equation model based on directional links indicated by the DAG. This approach, which relies on the seminal work by Pearl (2009) has been known to work well for the purpose of illuminating direct and indirect effects, the size of which can then be estimated once a hypothetical model is obtained from the DAG (see for example work by Bessler, 2001, 2002). This is particularly useful in cases, such as here, where theory is scant or non-existent as a guide as to which predictors drive new prescriptions. To give an idea of how this approach works, we use a modified PC (Partial Correlation) algorithm (PC Pattern) from the Tetrad IV package (The Tetrad Project 2014) to construct a DAG from the synthetic data, with the resulting graph displayed in Figure 4. We also give a brief description of the PC algorithm.

A few interesting features emerge from Figure 4. Samples seem to have an effect on new prescriptions, but indirectly, via their retail value. The cost per contact has an effect on the number of contacts, which seems sensible. We note that the PC algorithm typically yields not one but a class of directed acyclic graphs which are compatible with the data. When a directed arrow appears on the graph, this means that the direction was the same in all graphs generated by the algorithm. When this is not the case, the direction of the link is ambiguous and a link appears without arrows, as for example between *sam* and *eus*. The only variable which links directly into *nrx* is *rvs*, the retail value of samples.

*Figure 4: Directed Acyclic Graph for twelve synthetic variables*

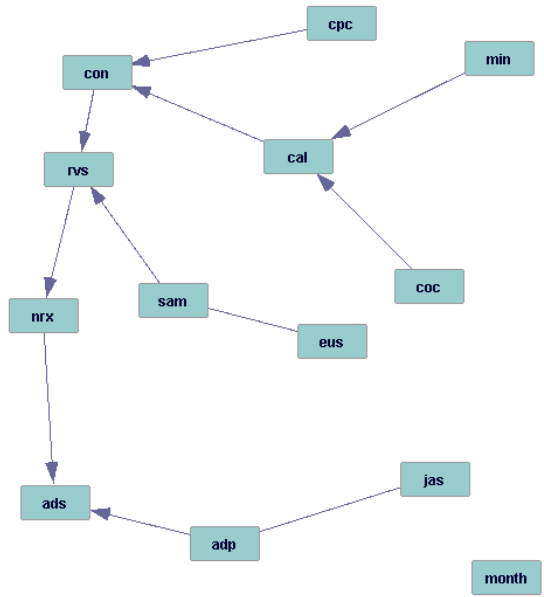

**Brief description of the PC (Partial Correlation) algorithm**

The PC algorithm works essentially as follows. First, a complete undirected graph is created with each variable corresponding to a vertex. Then edges are removed in pairs with variables which are independent, either unconditionally or conditionally on a subset of the remaining variables. Independence is tested with standard correlation tests (for continuous data assumed to be multivariate normal) or by using a test of independence in contingency tables (for categorical data). Note that Tetrad allows for either continuous data, or categorical data, but not a mixture of both types of data, unless the user provides a-priori information on which pairs of variables are independent conditionally on other variables. For a discussion of Tetrad and other acyclic graph software packages, see Haughton, Kamis and Scholten, 2006.

In order to orient the surviving links, the PC algorithm proceeds as follows: for each triplet *X*, *Y*, *Z* such that both pairs *X*, *Y* and *Y*, *Z* are linked but the pair *X*, *Z* is not linked, if *Y* does not appear in any set of variables which when conditioned upon makes *X* and *Z* independent, then the triplet *X*, *Y*, *Z* is oriented *X*→ *Y*← *Z*, in effect making *Y* a collider. This makes sense since such an orientation implies that *X* and *Z* are dependent given *Y*. Once all such colliders are identified, the algorithm proceeds like so: if *X*→ *Y*, *Y* and *Z* are linked and *X* and *Z* are not linked, and if there is no arrowhead at *Y* from *Z*, then *Y*, *Z* is oriented as *Y*→ *Z*. Such an orientation implies that *X* and *Z* are independent given *Y*. Appendix B of the Tetrad III user manual (The Tetrad Project 2010) explains how the algorithm unfolds on a particular example. See also Haughton and Haughton (2011), chapter 5, for an introduction to directed acyclic graphs with examples.

It is important to note that Directed Acyclic Graphs can infer causality as defined in Pearl (2009), but under very restrictive assumptions which tend not to be satisfied in many cases. So one would not want to claim any hope of getting definite causality statements out of this analysis; however what is obtained is a better understanding of which effects might be direct effects and which might be indirect.

**Discussion**

This paper has proposed a road map which should be helpful to researchers and practitioners involved in media mix modeling. The synthetic data on which our analyses rely relate to the marketing of a drug, but the methodology can be used in contexts other that health care as well.

The issue of negative coefficients in media mix models is one which, while unmentioned in the literature, is quite an important practical problem. In many real-life situations, negative coefficients are either ignored, or their sign ignored, when establishing the importance of each component. We have illustrated the fact that the key issue with negative coefficients is whether they arise because of multi-collinearity in the model or because of the genuine negative effect of a particular component. Discovering which of these two options holds is important.

Many media mix models include expenditures which are quite highly correlated. Because the mechanism which yields new prescriptions or higher sales in general is likely to be quite complex, we have suggested in this paper that Directed Acyclic Graphs can be pressed into service to help unravel which expenditures are directly linked to new prescriptions and which expenditures might act via an intermediary variable.

Several issues of interest remain for future work. One is the matter of the presence of competitor promotion expenditures when a few major competitors are present in a particular market. We suggest and investigate elsewhere (Haughton et al. 2013) that ignoring competing promotion activity in a media mix model could lead to biased results. Assuming that data on competing sales are available, it is possible to impute with good success the presence or absence of competitor promotional activity when two competitors are in place via the use of hidden Markov models (Haughton et al. 2013). More work remains to be done on situations with more than two main competitors, for example.

It also would be very interesting to apply these techniques to a real life data set. The synthetic data set presented in this paper is very representative of data used in a media mix marketing effort for the marketing for a drug, but nevertheless, applying the methods to actual data would be very valuable. The difficulty here is the scarcity of such datasets available to researchers.

In particular, a real life data set is very likely to include time series with auto-correlated variables. An interesting extension of this paper would be to adjust the methods to take into account auto-correlation in the errors of the media mix regression model.

Last but not least, the importance of word of mouth and viral marketing, notably via social networks, and in particular via social networks which are related to the medical profession, inclusive of co-publication networks, cannot be underestimated. The problem of how to combine social network data (on influential prescribers etc) with more traditional data used in media mix modeling is becoming quite pressing. Media mix models of the future will very probably need to incorporate social network components, since ignoring these components risks leading to biased results. A rather key challenge there is to link social network data to existing corporate databases used for targeting physicians and for building marketing programs.

**Acknowledgments**

The authors would like to thank the editor and two referees for most useful comments.